\documentclass{article}
\usepackage{epsfig}
\usepackage{amsthm}
\usepackage{amssymb}
\usepackage{xcolor}
\input mssymb.tex

\newtheorem{proposition}{Proposition}
\newtheorem{example}{Example}

\newtheorem{definition}{Definition}
\newtheorem{lemma}{Lemma}
\newtheorem{corollary}{Corollary}
\newtheorem{remark}{Remark}

\def\R{\Bbb R}

\def\spp{\vspace{5pt}  \noindent}

\def\dk{\partial K}

\usepackage{graphicx}
\begin{document}
\title{Convex Obstacles from Travelling Times} 


\author{Lyle Noakes\thanks{Department of Mathematics, University of Western Australia, Crawley WA 6009,
Australia (lyle.noakes@uwa.edu.au)}
\and Luchezar Stoyanov\thanks{Department of Mathematics, University of Western Australia, Crawley WA 6009,
Australia (luchezar.stoyanov@uwa.edu.au)}}

\date{}
\maketitle
\thispagestyle{empty}


{\bf Abstract:}  A construction is given for the recovery of a disjoint union of strictly convex smooth planar obstacles from travelling-time information. The obstacles are required to be such that no Euclidean line meets more than two of them.  

\spp


\spp


\spp
\section{Introduction}\label{sec0}

\spp
For some $n\geq 1$, let $K_1,K_2,\ldots ,K_n$ be disjoint closed convex subsets of Euclidean $2$-space 
$E^2\cong \mathbb{R}^2$, with each boundary $\partial K_k$ a $C^\infty$ strictly convex Jordan curve. 
Let  $K:=\cup_{k=1}^nK_i$ be contained in the interior of the bounded component $B$ of $E^2\setminus C$, 
where $C\subset E^2$ is also a strictly convex Jordan curve.
%
%

\spp
By a {\em geodesic} in the closure $M$ of $E^2\setminus K$ we mean a piecewise-affine constant-speed curve 
$x:\R \rightarrow M$ whose junctions are points of reflection on $\partial K$. The restriction of $x$ to an interval is also called a geodesic, and the set of all geodesics $x:[0,1] \rightarrow M$ with 
$x(0),x(1)\in C$ is denoted by ${\mathcal X}$. Then $x\in {\mathcal X}$ is critical for the {\em length functional} 
$$J(x) ~:=~\int _0^1\Vert \dot x(t)\Vert ~dt$$
over constant-speed piecewise-$C^1$ curves in $M$ satisfying $x(0)=x_0,x(1)=x_1$.
By Theorem 1.1 of \cite{NS1}, $K$ is uniquely determined by its {\em travelling-time data}
$${\mathcal T}~:=~\{ (x(0),x(1),J(x)):x\in {\mathcal X}\} .$$
Similar results are proved in \cite{NS2} for obstacles in $E^m$ where $m>2$. Unfortunately the proof of Theorem 1.1 in \cite{NS1} is not constructive: all that is shown is that ${\mathcal T}$ is different for 
different convex obstacles $K$. When $n=1$ it is straightforward to calculate $K$ from ${\mathcal T}$, and with a little more effort $K$ can also be reconstructed when $n=2$ (see Sect. 4 in \cite{NS3}). 
More  interestingly\footnote{Indeed \cite{St2}  applies in a much more general setting, where the obstacles are not necessarily convex, and $E^2$ is replaced by a Riemannian manifold of any finite 
dimension.}, Theorem 1.1 of \cite{St2}  allows to compute the  area of $K$  from ${\mathcal T}$. Importantly, the application of the result in \cite{St2} is made possible by the fact 
(proved in \cite{NS1}) that the set of  points generating trapped trajectories in the exterior of obstacles
$K$ considered in this paper has  Lebesgue measure zero. Constructing $K$ is equivalent to constructing $\partial K$, but this seems difficult for $n\geq 3$. 
In the present paper we show how to construct $\partial K$ from ${\mathcal T}$ when $K$ is  {\em in general position}, namely when no line meets more than two connected components of $K$. 
 
 \spp
 Inverse problems concerning metric rigidity have been studied for a long time in Riemannian geometry: we refer to  
\cite{SUV}, \cite{Gu} and their references for more information. In the last 20 years or so similar problems have been considered for scattering by obstacles, where the task is to recover geometric information about an obstacle from its scattering length spectrum \cite{St1}, or from travelling times of
scattering rays in its exterior \cite{NS3}. 

\spp
In general, an obstacle in the Euclidean space $E^m\cong \mathbb{R}^m$ ($m\geq 2$) is a compact subset $K$ of ${\R}^m$ 
with a smooth (e.g. $C^{3}$) boundary $\partial K$ such that  $\Omega_K = \overline{{\R}^m\setminus K}$ is connected. 
The { scattering rays} in  $\Omega_K$ are generalized geodesics  (in the sense of Melrose and  Sj\"ostrand \cite{MS1}, \cite{MS2}) that are unbounded in both directions.  
Most of these scattering rays are billiard trajectories  with finitely many reflection points at $\dk$. When $K$ is a finite disjoint union of  strictly convex domains, then all scattering rays in 
$\Omega_K$ are billiard trajectories, namely geodesics of the type described above.

\spp
It turns out that some kinds of obstacles are uniquely recoverable from their travelling times spectra. For example, as mentioned above, this was proved in \cite{NS2} for obstacles  $K$ 
in $\R^m$ ($m \geq 3$) that are finite disjoint unions of strictly convex bodies with $C^3$ boundaries. The case $m = 2$ requires a different proof, given recently in \cite{NS1}.

\spp
The set of the so called trapped points (points that generate trajectories with infinitely many reflections)  plays a rather important role in various inverse problems in scattering by 
obstacles, and also in problems on metric rigidity in Riemannian geometry.  As an example of M. Livshits shows   (see e.g. Figure 1  in \cite{NS4} or \cite{St2}),  in general the set of trapped points
may contain a non-trivial open set. In such a case the  obstacle cannot be recovered from travelling times. In dimensions $m > 2$ examples similar to that of Livshits were given  in \cite{NS4}.

 \spp
 The layout of the paper  is as folllows. 
 
 \spp
 In \S \ref{sec1} we collect some simple observations about linear (non-reflected) geodesics. This leads to the construction of $4(n^2-n)$ so-called {\em vacuous arcs} $\beta _j$ in ${\mathcal T}$, and then 
 $2(n^2-n)$ {\em initial arcs} in $\partial K$. Our plan is to build on the initial arcs, using travelling-time data from reflected rays to construct {\em incremental arcs} in $\partial K$, until 
 eventually\footnote{Unlike the  initial arcs, there are countably many incremental arcs, yielding diminishing additional information from ever-increasing amounts of precisely known data. In practice, 
 insufficient data and limited computing power makes  it difficult to carry out more than a few inductive steps, and $\partial K$ is found only approximately.} the whole of $\partial K$ is found. In \S \ref{prsec} 
 we describe an {\em inductive step} for constructing incremental arcs from previously determined arcs, and from observations of ${\mathcal T}$. To make the relevant observations we need to understand 
 some of the mathematical structure of ${\mathcal T}$. 
 
 \spp
The first step towards this understanding is made in \S \ref{nonvac}, where some simple facts about (typically non-reflected) geodesics are recalled. These facts, including a known result for computing initial directions 
of geodesics, are applied in \S \ref{nowtgt} to investigate the structure of travelling-time data of nowhere-tangent geodesics. In particular, cusps in so-called {\em telegraphs} of ${\mathcal T}$ correspond to geodesics 
that are tangent to $\partial K$. 

\spp
The family of all such cusps is studied in \S \ref{singtgt}, where the augmented travelling-time data $\tilde {\mathcal T}$ is shown to be the closure of a countable family of disjoint open $C^\infty$ arcs $\tilde \beta _j$. As described 
in \S \ref{prsec}, the property of {\em extendibility} can be checked for each  $\tilde \beta _j$. When $\tilde \beta _j$ is extendible it yields an incremental arc in $\partial K$. When $\tilde \beta _j$ is not extendible, a 
trick using general position replaces $\tilde \beta _j$ by an extendible $\tilde \beta _{\bar j}$ yielding an incremental arc as previously.

\section{Linear Geodesics and Vacuous Arcs}\label{sec1}

\spp
From now on let $K = \cup_{k=1}^nK_i$ be an obstacle in $E^2$,  where $K_1,K_2,\ldots ,K_n$ are disjoint closed convex subsets of $E^2$
with boundaries that are $C^\infty$ strictly convex Jordan curves. As before, assume that $K$ is contained in the interior of the bounded component $B$ 
of $E^2\setminus C$,  where $C\subset E^2$ is also a strictly convex Jordan curve. We also assume that $K$ is in general position.

\begin{lemma}\label{lem0} Geodesics in $\mathcal{X}$ are not tangent to $\partial K$, except perhaps at the first or 
last points of contact with $\partial K$ (either or both). 
\end{lemma}

\spp
{\bf Proof:} If tangency was at an intermediate point of contact, the tangent line would have common points with at least $3$ connected components of $K$, contradicting general position. 
\qedsymbol

\spp
We begin by investigating travelling times of {\em linear geodesics}, namely geodesics in ${\mathcal X}$ that do not reflect at all. The travelling time data from linear geodesics is  
$${\mathcal T}_0:={\mathcal T}\cap \{ (x_0,x_1,\Vert x_1-x_0\Vert ):x_0,x_1\in C\} $$ 
from which we find 
$$\tilde \partial {\mathcal{T}_0}:=(C\times C\times (0,\infty ))\cap \partial {\mathcal T}_0 =\cup_{q\geq 1}\mathcal{T}_0^q$$
where ${\mathcal T}_0^q$ is defined as  the travelling-time data from geodesics meeting $\partial K$ exactly $q$ times tangentially and nowhere else.  
By Lemma \ref{lem0} $\mathcal{T}_0^q=\emptyset$ for $q\geq 3$, namely~   
$\displaystyle{\tilde \partial {\mathcal T}_0={\mathcal T}_0^1\cup {\mathcal T}_0^2}$. In the simplest case where $n=1$, $\mathcal{T}_0^2$ is empty and $\partial K$ is constructed as  the envelope of the line segments $[x_0,x_1]$ 
where $(x_0,x_1,\Vert x_1-x_0\Vert )\in \tilde \partial \mathcal{T}_0$. Suppose $n\geq 2$ from now on. 

\begin{proposition}\label{prop1} ${\mathcal T}_0^1$ is a union of $4(n^2-n)$ nonintersecting  bounded open $C^\infty$ arcs $\beta _j$ whose boundaries in $\tilde \partial {\mathcal T}_0$ comprise 
${\mathcal T}_0^2$ which is finite of size $4(n^2-n)$. 
\end{proposition}

\spp
{\bf Proof:} For $1\leq k\not= k'\leq n$ there are $8$ {\em directed} Euclidean line segments (linear bitangents) tangent to both $\partial K_k$ and $\partial K_{k'}$. Each directed linear bitangent is an endpoint of two 
maximal open arcs of  directed line segments that are singly-tangent. The travelling-time data for the linear bitangents is $\mathcal{T}_0^2$. The 
 travelling-time data for the open arcs $\beta _j$ are the path components of $\mathcal{T}_0^1$. 
\qedsymbol

\spp
So $n$ is found from $\mathcal{T}_0^1$. 
\begin{definition} For $1\leq j\leq 4(n^2-n)$ the {\em conjugate} $\bar j$ is defined to be $j+2(n^2-n)$ or $j-2(n^2-n)$ according as $1\leq j\leq 2(n^2-n)$ or $2(n^2-n)+1\leq j\leq 4(n^2-n)$. \qedsymbol
\end{definition}

\spp
Evidently $\bar{\bar j}=j$. Order the arcs $\beta _j$ in ${\mathcal T}_0^1$ so that 
$(x_0,x_1,t)\in \beta _j~\Leftrightarrow (x_1,x_0,t)\in \beta _{\bar j}$. The {\em initial arcs} in $\partial K$ are the $2(n^2-n)$ nonempty disjoint connected open subsets of $\partial K$ found as the envelopes of the $[x_0,x_1]$, 
where $(x_0,x_1,t)\in \beta _j$ for $1\leq j\leq 2(n^2-n)$. For each $1\leq k\leq n$ there are $2(n-1)$ initial arcs in $\partial K_k$.

\section{Nonlinear Geodesics}\label{nonvac}

\spp
In order to construct envelopes of other singly-tangential geodesics, we shall identify the travelling-time data ${\cal T}^q$ of $q$-times tangential geodesics in ${\mathcal X}$, especially $q=1$. Whereas  
${\mathcal T}_0^1$ is found by simple inspection of ${\mathcal T}$,  some effort is required to isolate ${\cal T}^1$. We first recall some known results about directions of geodesics and travelling-times. 

\spp  
For $(x_0,v_0)\in (E^2-K)\times E^2$ let $x_{x_0,v_0}:\R \rightarrow M$ be the geodesic satisfying $x_{x_0,v_0}(0)=x_0$ and $\dot x_{x_0,v_0}(0)=v_0$. The {\em endpoint map} 
${\mathcal E}:(E^2-K)\times E^2 \rightarrow M$ is the continuous function given by ${\mathcal E}(x_0,v_0):=x_{x_0,v_0}(1)$. 
\begin{lemma}\label{lem1}
For $x_0\in E^2-K$,  suppose that $x_{x_0,v_0}\vert [0,1]$  is nowhere tangent to $\partial K$, and that $x_{x_0,v_0}(1)\notin \partial K$. Then ${\mathcal E}$ is $C^\infty$ near $(x_0,v_0)\in (E^2-K)\times E^2$, 
and the restriction of its derivative  $d{\mathcal E}_{x_0,v_0}:\R ^2\times \R ^2 \rightarrow \R^2$ to $\{ {\bf 0}\} \times \R ^2$ is a linear isomorphism. 
\end{lemma}

\spp
{\bf Proof:} Because the $K_k$ are disjoint and strictly convex, the endpoints of $x_{x_0,v_0}\vert [0,1]$ are nonconjugate. 
\qedsymbol

\spp 
\begin{lemma}\label{lem2} For $x_0\not= x_1:={\mathcal E}(x_0,v_0)$ and with the hypotheses of Lemma \ref{lem1}, there exists an open neighbourhood $U_0$ of $x_0$ in $E^2-K$, and a unique $C^\infty$ function $\phi :=\phi _{x_0,v_0}:U_0\rightarrow (0,\infty )$ whose gradient ${\nabla} \phi$ is everywhere of unit length,  satisfying  
$${\mathcal E}(\tilde x_0,-\phi (\tilde x_0){\nabla} \phi (\tilde x_0))=x_1$$
for all $\tilde x_0\in U_0$. Here $\phi (x_0)=\Vert v_0\Vert$ with ${\nabla}\phi (x_0)=-v_0/\vert v_0\Vert$.
 \end{lemma}

\spp
{\bf Proof:} By Lemma \ref{lem1} and the implicit function theorem, there exist a unique $C^\infty$ function $X:U_0\rightarrow E^2$ satisfying 
${\mathcal E}(\tilde x_0,X (\tilde x_0))=x_1$ for all $\tilde x_0\in U_0$. 
Because $x_0\not= x_1$, $X$ is never-zero for  $U_0$ sufficiently small.  
Then the geodesic $x_{\tilde x_0,X(\tilde x_0)}:[0,1]\rightarrow M$ joining $\tilde x_0,x_1$, has length 
$$\displaystyle{\phi (\tilde x_0):=\Vert X(\tilde x_0)\Vert =J(x_{\tilde x_0,X(\tilde x_0)})}.$$ 
Differentiating with respect to $\tilde x_0\in U_0$ in the direction of  $\delta \in \R ^2$, we find  
$d\phi _{\tilde x_0}(\delta) =-\langle X(\tilde x_0)/\Vert X(\tilde x_0)\Vert ,\delta \rangle$, because geodesics are critical for $J$ when variations have {\em fixed} endpoints.   
\qedsymbol

\spp
The {\em order} $o(x)$ of a geodesic $x$ is the number of intersections with $\partial K$. Write ${\mathcal X}_{r}:=\{ x\in {\mathcal X}:o(x)=r\}$. 

\spp
Let $d_K$ be the minimum distance between obstacles. 
Writing $\tau (x):=t$ for the travelling time (length) of $x\in {\mathcal X}$,
\begin{equation}\label{bds0}o(x)d_K\leq J (x)\leq o (x){\rm diam}(C).\end{equation} 
\section{Arcs and Generators for Nowhere-Tangent Geodesics}\label{nowtgt}

\spp
Let ${\mathcal X}^q$ be the space of geodesics $x\in {\mathcal X}$ that  are exactly $q$-times tangent to $\partial K$, and set 
$${\mathcal T}^q:=\{ (x(0),x(1),J(x)):x\in {\mathcal X}^q\} \subset {\mathcal T}.$$
Then ${\mathcal T}_0^q\subset {\mathcal T}^q$ and, by Lemma \ref{lem0},  ${\mathcal T}^q=\emptyset$ for $q\geq 3$. We have constructed $ {\mathcal T}_0^q$ from $\mathcal{T}$,  but not yet ${\mathcal T}^q$ for $q\leq 2$. 
For $x_1\in C$ define 
$${\mathcal T}_{x_1}:=\{ (x_0,t):(x_0,x_1,t)\in {\mathcal T}\} \subset C\times \R ,~~{\mathcal T}_{x_1}^q:=\{ (x_0,t): (x_0,x_1,t)\in {\cal T}^q\} \subset {\mathcal T}_{x_1} .$$
Likewise $\mathcal {X}_{x_1}\subset \mathcal{X}$ is the set of geodesics $x:[0,1]\rightarrow M$ with $x(1)=x_1$. Write ${\mathcal X}_{x_1}^q:={\mathcal X}_{x_1}\cap {\mathcal X}^q$. 
\begin{remark} ${\mathcal T}_{x_1}$ is found directly from ${\mathcal T}$, but ${\mathcal T}_{x_1}^q$ is yet to be determined for $q\leq 2$. \qedsymbol
\end{remark}
\begin{remark} ${\mathcal T}_{x_1}^0$ is open and dense in ${\mathcal T}_{x_1}$, ${\mathcal T}_{x_1}^1$ is open and dense in ${\mathcal T}_{x_1}^1\cup {\mathcal T}_{x_1}^2$, and ${\mathcal T}_{x_1}^2$ is discrete. \qedsymbol
\end{remark}
\begin{proposition}\label{prop11} For $x_1\in C$  we have 
 ${\mathcal T}_{x_1}^0=\cup _{i\geq 1}\alpha _{i,x_1}$, where
 \begin{enumerate}
 \item the $\alpha _i:=\alpha _{i,x_1}$ are pairwise-transversal $C^\infty$ open bounded arcs in 
$C\times \R$, 
\item  $\cup _{q\geq 1}{\mathcal T}_{x_1}^q=\cup _{i\geq 1}\partial \alpha _i$,
\item each $\alpha _i$ has a {\em generator}, namely a $C^\infty$ function $\phi _{i}:U_i\rightarrow \R$ with $U_i\subset E^2$ open, such that
\begin{itemize}
\item $U_i\cap C$ is an open arc in $C$, and $\tilde x_0\mapsto (\tilde x_0,\phi _{i}(\tilde x_0))$ is a diffeomorphism from $U_i\cap C$ onto $\alpha _i$,
\item $x_{\tilde x_0,\nu _i (x_0)}\vert [0,1]\in {\cal X}_{x_1}^0$, where $\nu _i(\tilde x_0):=-\phi _i(\tilde x_0)\nabla \phi _{i} (\tilde x_0)$ for $\tilde x_0\in U_i\cap C$.
\end{itemize}
\end{enumerate}
\end{proposition} 

\spp
{\bf Proof:} For $(x_0,t)\in {\mathcal T}_{x_1}^0$ there exists $x\in {\mathcal X}_{x_1}$ with $J(x)=t$.  
Then $x_1={\mathcal E}(x_0,v_0)$ 
where $v_0=\dot x(0)$. By Lemma \ref{lem2}, for some open neighbourhood $U_0$ of $x_0$ in $E^2$, 
there is a unique $C^\infty$ function $\phi =\phi _{x_0,v_0}:U_0\rightarrow (0,\infty )$ with $\phi (x_0)=t$ and $\nabla \phi _{x_0}=-v_0/\Vert v_0\Vert$, such that  
${\mathcal E}(\tilde x_0,-\phi (\tilde x_0)\nabla \phi (\tilde x_0))=x_1$ for all $\tilde x_0\in U_0$. In particular the last equation holds for $\tilde x_0\in U_0\cap C$, namely $\tilde x_0\mapsto (\tilde x_0,\phi (\tilde x_0))$ 
embeds $U_0\cap C$ in ${\mathcal T}_{x_1}^0$. By continuation, the embedding extends uniquely in both directions around $C$, until just before $x$ is tangent to some $\partial K_k$, which must eventually happen. 
So ${\mathcal T}_{x_1}^0$ is a countable union of $C^\infty$ embedded arcs $\alpha _i$. 
 
\spp
Pairwise transversality is proved by contradiction as follows. Suppose $\alpha _i\not= \alpha _{i'}$ meet tangentially at $(x_0,t)\in {\mathcal T}_{x_1}^0$. Then 
$\langle \nabla \phi _{x_0,v_{i}}(x_0),w\rangle =\langle \nabla \phi _{x_0,v_{i'}}(x_0),w\rangle $ where $\Vert v_i\Vert =t=\Vert v_{i'}\Vert $, and $w\not= {\bf 0}$ is tangent to $C$ at $x_0$. By Lemma \ref{lem2}, 
$\nabla \phi _{x_0,v_{i}}(x_0)$ and $\nabla \phi _{x_0,v_{i'}}(x_0)$ point out from $B$ at $x_0$. So $-v_i=t\nabla \phi _{x_0,v_{i}}(x_0)=t\nabla \phi _{x_0,v_{i'}}(x_0)=-v_{i'}$ by Lemma \ref{lem2}, contradicting 
$\alpha _i\not=\alpha _{i'}$.  

\spp
Because ${\mathcal T}_{x_1}^0$ is open and dense in ${\mathcal T}_{x_1}$,  $\cup _{q\geq 1}{\mathcal T}_{x_1}^q=\cup _{i\geq 1}\partial \alpha _i$. \qedsymbol

\spp
By continuity, the 
orders $o(\alpha _i)$ of the $x_{\tilde x_0,\nu  _{i} (\tilde x_0)}\in {\mathcal X}^0$ 
are independent of $\tilde x_0\in U_i\cap C$. From (\ref{bds0}) we obtain, for all $\tilde x_0\in U_i\cap C$,   
\begin{equation}\label{bds}o(\alpha _i)d_K\leq \phi _{i}(\tilde x_0)\leq o (\alpha _i){\rm diam}(C),\end{equation} 
and the arcs $\alpha _i$ are similarly bounded. 
For any $i$, the closures $\bar \alpha _i$ and $\bar \alpha _{i'}$ in ${\mathcal T}_{x_1}$ are disjoint for all but finitely many $i'$, where $i,i'\geq 1$. 
The generator $\phi _i$ defines $x_{\tilde x_0,\nu _i(\tilde x_0)}\in {\mathcal X}_{x_1}^0$ for every $\tilde x_0\in U_i\cap C$. For $(x_0,t)\in \partial \alpha _i$, define 
$$\nu _i(x_0):={\rm lim}_{\tilde x_0\rightarrow x_0}\nu _i(\tilde x_0)\in E^2.$$
Then $x_{x_0,\nu _i(x_0)}\in {\mathcal X}_{x_1}^1\cup {\mathcal X}_{x_1}^2$, and $(x_0,t )\in {\cal T}_{x_1}^1\cup {\cal T}_{x_1}^2$ where $t=\Vert \nu _i(x_0)\Vert $. 

\spp
\begin{proposition}\label{cor2} If $(x_0,t)\in {\mathcal T}_{x_1}^1$ 
then $\{ (x_0,t)\} =(\partial \alpha _i )\cap (\partial \alpha _{i'})$ for some unique $i, i'$ with $o(\alpha _{i'})=o(\alpha _i)+1$.   Then $U_i\cap C$ and $U_{i'}\cap C$ on the same side of $x_0$ in $C$ and, 
for $\tilde x_0\in U_i\cap U_{i'}\cap C$, 
$\phi _{i'}(\tilde x_0)>\phi _{i}(\tilde x_0)$.  We also have  
$$\lim_{\tilde x_0\rightarrow x_0}\phi _i(\tilde x_0)=\lim_{\tilde x_0\rightarrow x_0}\phi _{i'}(\tilde x_0)=t\hbox{ ~and~ } 
\lim_{\tilde x_0\rightarrow x_0}\nabla \phi _i(\tilde x_0)=\lim_{\tilde x_0\rightarrow x_0}\nabla \phi _{i'}(\tilde x_0).$$
\end{proposition}

\spp
{\bf Proof:}  We can write $t=J(x_{x_0,v_0})$ where $x_{x_0,v_0}\in {\mathcal X}_{x_1}^1$ and $\Vert v_0\Vert =t$. 
Suppose the last (respectively first) segment of $x_{x_0,v_0}$ is not tangent to $\partial K$. Then, by general position, 
the first (last) segment is tangent. Perturbing the last segment while maintaining the endpoint $x_1$, gives two arcs of nowhere-tangent geodesics, whose initial points $\tilde x_0$ lie on the same side of $x_0$ in $C$. 
Along one arc the first (last) segment remains linear and the order decreases by $1$. Along the other arc, the first (last) segment breaks into two linear segments, maintaining the order and increasing the travelling time. 

\spp
For $\tilde x_0$ near $x_0$, the two arcs of geodesics define arcs $\{ (\tilde x_0,\phi _i(\tilde x_0))\}$, $\{ (\tilde x_0,\phi _{i'}(\tilde x_0))\}$ in ${\mathcal T}_{x_1}^0$ contained in maximal arcs $\alpha _i,\alpha _{i'}$, labelled so that 
$o(\alpha _{i'})=o(\alpha _i)+ 1$. Then  $\nu _{i'}(x_0)=v_0=\nu_i(x_0)$, and   
$\phi _{i'}(\tilde x_0)>\phi _{i}(\tilde x_0)$. 
We also have $\lim_{\tilde x_0\rightarrow x_0}\phi _i(\tilde x_0)=\Vert v_0\Vert =\lim_{\tilde x_0\rightarrow x_0}\phi _{i'}(\tilde x_0)$, and~  
$\displaystyle{\lim_{\tilde x_0\rightarrow x_0}\nabla \phi _i(\tilde x_0)=v_0/\Vert v_0\Vert =\lim_{\tilde x_0\rightarrow x_0}\nabla \phi _{i'}(\tilde x_0)}$. ~\qedsymbol

\spp
Since ${\mathcal T}_{x_1}^1$ is dense in ${\mathcal T}_{x_1}^1\cup {\mathcal T}_{x_1}^2$, Proposition \ref{cor2} has the 
\begin{corollary} ${\mathcal T}_{x_1}^1\cup {\mathcal T}_{x_1}^2$ is the closure in ${\mathcal T}_{x_1}$ of the set of all points  $(x_0,t)\in {\mathcal T}_{x_1}$ where ${\mathcal T}_{x_1}$ has an isolated cusp. \qedsymbol
\end{corollary} 

\begin{figure}[h] 
   \centering
   \includegraphics[width=5in]{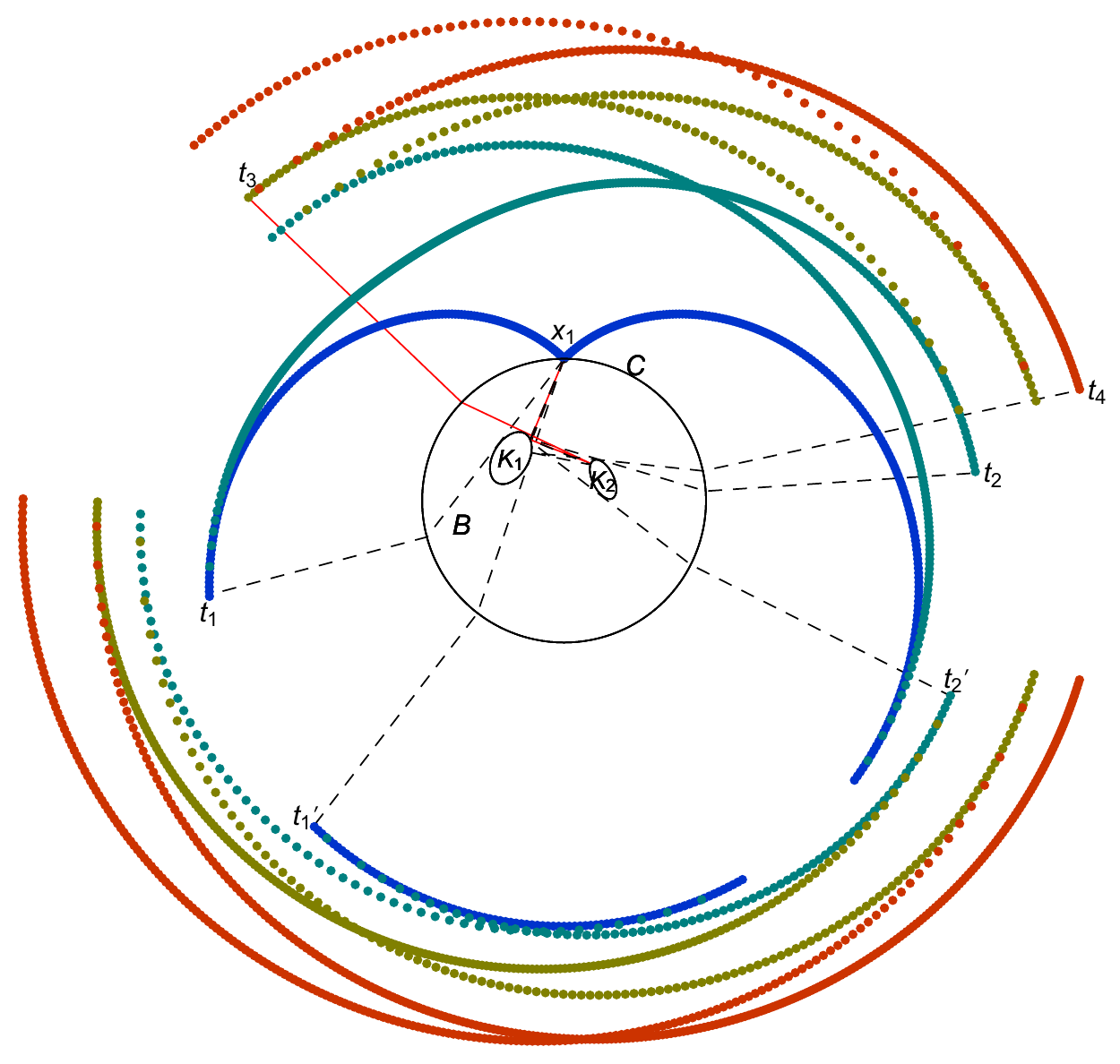} 
   \caption{{\rm Part of  the telegraph in Example \ref{ex1}}}
   \label{pic1}
\end{figure}

\spp
A $C^\infty$ embedding
$\epsilon$ of ${\mathcal T}_{x_1}$ in $E ^2$ is given by  
$\epsilon (x_0,t):=x_0+t\nu (x_0)$, 
with $\nu :C\rightarrow E^2$ some constant-length nonzero outward-pointing normal field.
Cusps in ${\mathcal T}_{x_1}$ are found by inspecting the {\em telegraph} at $x_1$, defined as $\epsilon ({\mathcal T}_{x_1})\subset E^2$. 
\begin{example}\label{ex1}
Figure \ref{pic1} displays part of $\epsilon ({\mathcal T}_{x_1})$ with $x_1=(0.4,4)$ with $n=2$, and $C$ the circle of radius $4$ and centre $(0.4,0)$. The telegraph is mainly smooth, but different arcs (light-blue, yellow-green and red) meet in cusps, and $6$ transversal self-intersections are seen. Cusps (labelled $t_1,t_1',t_2,t_2',t_3,t_4$) correspond to tangencies of geodesics ending at $x_1$ to $K_1$ or $K_2$.   
 \qedsymbol
\end{example}

 \spp
Next we augment ${\mathcal T}_{x_1}$ and ${\mathcal T}$ to data sets $\tilde {\mathcal T}_{x_1}$ and $\tilde {\mathcal T}$ that include initial velocities of geodesics. We first exclude points of intersection of the open arcs $\alpha _i=\alpha _{i,x_1}$ in Proposition \ref{prop11} (these points are reinserted later), by defining $\alpha _i^*=\alpha _{i,x_1}^*:=\alpha _{i}-\cup_{i'\not= i}\alpha _{i'}$. 
\begin{remark} Any $\alpha _i$ intersects at most finitely many $\alpha _{i'}$. Because intersections of $\alpha _i$ and $\alpha _{i'}$ are transversal for $i\not= i'$, ${\mathcal T}_{x_1}^{0*}:=\cup _{i\in I}\alpha _{i,x_1}^*$ is dense in  
${\mathcal T}_{x_1}^{0}:=\cup _{i\in I}\alpha _{i,x_1}$, 
and 
$${\mathcal T}^{0*}~:=~\{ (x_0,x_1,t):(x_0,t)\in {\mathcal T}_{x_1}^{0^*},~x_1\in C\} \hbox{~ is dense in ~}$$ 
$${\mathcal T}^{0}~:=~\{ (x_0,x_1,t):(x_0,t)\in {\mathcal T}_{x_1}^{0},~x_1\in C\} .$$
\qedsymbol
\end{remark}
\begin{remark} The $\alpha _{i,x_1}^*$ partition ${\mathcal T}_{x_1}^{0*}$. The generators $\phi _i$ restrict to $C^\infty$ functions  on the open subsets  
$$D_{i,x_1}:=\{ c\in C: (c,t)\in \alpha _{i,x_1}^*\}$$
of $C$. \qedsymbol
\end{remark}

\spp 
For $(x_0,t)\in {\mathcal T}_{x_1}^{0*}$ define $u_0=u_{x_0,t,x_1}$ to be the unit
vector $-\nabla \phi _i(x_0)$ pointing inwards from $C$. Then set
$$\tilde {\mathcal T}_{x_1}^{0*}~:=~\{ (x_0,u_0,t): (x_0,t)\in {\mathcal T}_{x_1}^{0*}\} \hbox{~ and ~}\tilde {\mathcal T}^{0*}~:=~\{ (x_0,u_0,x_1,t): (x_0,u_0,t)\in \tilde {\mathcal T}_{x_1}^{0*}\} $$
To reinsert the excluded points, define $\tilde {\mathcal T}^0$ to be the closure of $\tilde {\mathcal T}^{0*}$ in 
$$\displaystyle{\{ (x_0,u_0,x_1,t):(x_0,x_1,t)\in {\mathcal T}^0 \hbox{ with }u_0\in S^1\} },$$
and $\tilde {\mathcal T}_{x_1}^0:=  \{ (x_0,u_0,t):(x_0,u_0,,x_1,t)\in \tilde {\mathcal T}^0\}$. ~ Define 
$\tilde {\mathcal T}$ be the closure of $\tilde {\mathcal T}^0$ in 
$$\{ (x_0,u_0,x_1,t):(x_0,x_1,t)\in {\mathcal T} \hbox{ with }u_0\in S^1\} ,$$
and  $\tilde {\mathcal T}_{x_1}:=  \{ (x_0,u_0,t):(x_0,u_0,x_1,t)\in \tilde {\mathcal T} \}$. For $q\geq 1$ define 
$$\tilde {\mathcal{T}}^q:=
\{ (x_0,u_0,x_1,t)\in \tilde {\mathcal{T}}:(x_0,x_1,t)\in \mathcal{T}^q\}.$$

\section{Singly-Tangent Geodesics}\label{singtgt}

\spp
Summarising so far, for $x_1\in C$, 
\begin{itemize}
\item ${\mathcal T}_{x_1}$ is read directly from ${\mathcal T}$,
\item ${\mathcal T}_{x_1}^+:={\mathcal T}_{x_1}^1\cup {\mathcal T}_{x_1}^2$ (respectively ${\mathcal T}_{x_1}^0$) is the non-smooth (respectively smooth) part of  ${\mathcal T}_{x_1}$, 
\item we have seen how to find arcs $\alpha _{i,x_1}$ and generators $\phi _i$ for ${\mathcal T}_{x_1}^0$, 
\item $\tilde {\mathcal T}_{x_1}^+:=\tilde {\mathcal T}_{x_1}^1\cup \tilde {\mathcal T}_{x_1}^2$ and $\tilde {\mathcal T}_{x_1}^0$ are obtained using the $\phi _i$,  
\item $\tilde {\mathcal T}^+:=\tilde {\mathcal T}^1\cup \tilde {\mathcal T}^2$ and $\tilde {\mathcal T}^0$ are found by varying $x_1$.
\end{itemize}

\spp 
To distinguish 
 $\tilde {\mathcal T}^1$ from $\tilde {\mathcal T}^2$ we need Proposition \ref{prop4}, which is a structural result, analogous to Proposition \ref{prop11}. A geodesic $x^*\in {\mathcal X}$ is said to be {\em bitangent} when it has two points of tangency to $\partial K$. We call $x^*$ {\em linear} when it has no other points of contact with $\partial K$. 

\begin{proposition}\label{prop4} For some countable locally finite family $\tilde {\mathcal B}=\{ \tilde \beta _j:j\geq 1\}$ of disjoint bounded open 
$C^\infty$ arcs in $\tilde {\mathcal{T}}^+$,  
\begin{enumerate}
\item $\tilde {\mathcal T}^{1}=\cup_{j\geq }\tilde \beta _j$, 
\item for $1\leq j\leq 4(n^2-n)$, 
$\tilde \beta _j = \{ (x_0,(x_1-x_0)/t,x_1,t): (x_0,x_1,t)\in \beta _j\}$, where the $\beta _j$ are the vacuous arcs in ${\mathcal T}_0^1$, defined in \S \ref{sec1},~
\item for every\footnote{Including possibly $j>4(n^2-n)$.} $j\geq 1$ there is a diffeomorphism $\psi _j:V_j\rightarrow \tilde \beta _j$ where $V_j$ is an open arc in $C$, and $\psi _j(x_0)\in \{ x_0\} \times S^1\times C\times (0,\infty )$ for all $x_0\in V_j$,
\item each $(x_0^*,u_0^*,x_1^*,t^*)\in \tilde {\mathcal T}^2$ is 
an endpoint of four open arcs $\tilde \beta _j,\tilde \beta _{j'},\tilde \beta _{j''},\tilde \beta _{j'''}$, where three of $V_j,V_{j'},V_{j''},V_{j'''}$ are on one side of $x_0^*\in C$, and one is on the other side. 
\item $\tilde {\mathcal T}^2=\cup _{j\geq 1}\partial \tilde \beta _j$. 
\end{enumerate} 
\end{proposition}

\spp
{\bf Proof:} For $(x_0,u_0,x_1,t)\in \tilde {\mathcal{T}}^1$, we have $x_{x_0,tu_0}\in \mathcal{X}^1$ and $x_{x_0,tu_0}(1)=x_1$. 
Now $x_{x_0,tu_0}$ is tangent to $\partial K$ at precisely one point. By Lemma \ref{lem0} this is either the first or last point of contact with $\partial K$.  

\spp
If the tangency is first then, perturbing the point of tangency in $\partial K$ gives a small open $C^\infty$ arc around $(x_0,tu_0,x_1,t)$ contained in $\tilde {\mathcal{T}}^1$. Similarly, if the tangency is last, an open $C^\infty$ arc in $\tilde {\mathcal T}^+$ is given by perturbing the point of tangency in $\partial K$. 
So the path components $\tilde \beta _j$ of $\tilde {\mathcal{T}}^1$ in $\tilde {\mathcal T}^+$ are connected smooth $1$-dimensional submanifolds of $C\times S^1\times C\times \R$. They are bounded, nonclosed and,  
for $1\leq j\leq 4(n^2-n)$, can be listed as augmentations of the $\beta _j$. Then  {\it 1.} and {\it 2.} hold. 

\spp
For $(x_0^*,u_0^*,x_1^*,t^*)\in \tilde {\mathcal{T}}^2$ the geodesic $x^*=x_{x_0^*,t^*u_0^*}$ is tangent to $\partial K$ at both first and last points of contact, and nowhere else. Nearby geodesics in $\mathcal{X}^1$ are obtained by maintaining tangency either at a variable first point of contact, or at a variable last point of contact with $\partial K$. The tangencies at first (respectively last) points of contact generate arcs $\tilde \beta _j,\tilde \beta _{j'}$ (respectively $\tilde \beta _{j''},\tilde \beta _{j'''}$) in $\tilde {\mathcal T}^1$, separated by $(x_0^*,u_0^*,x_1^*,t^*)$.  

\begin{figure}[h] 
   \centering
   \includegraphics[width=3in]{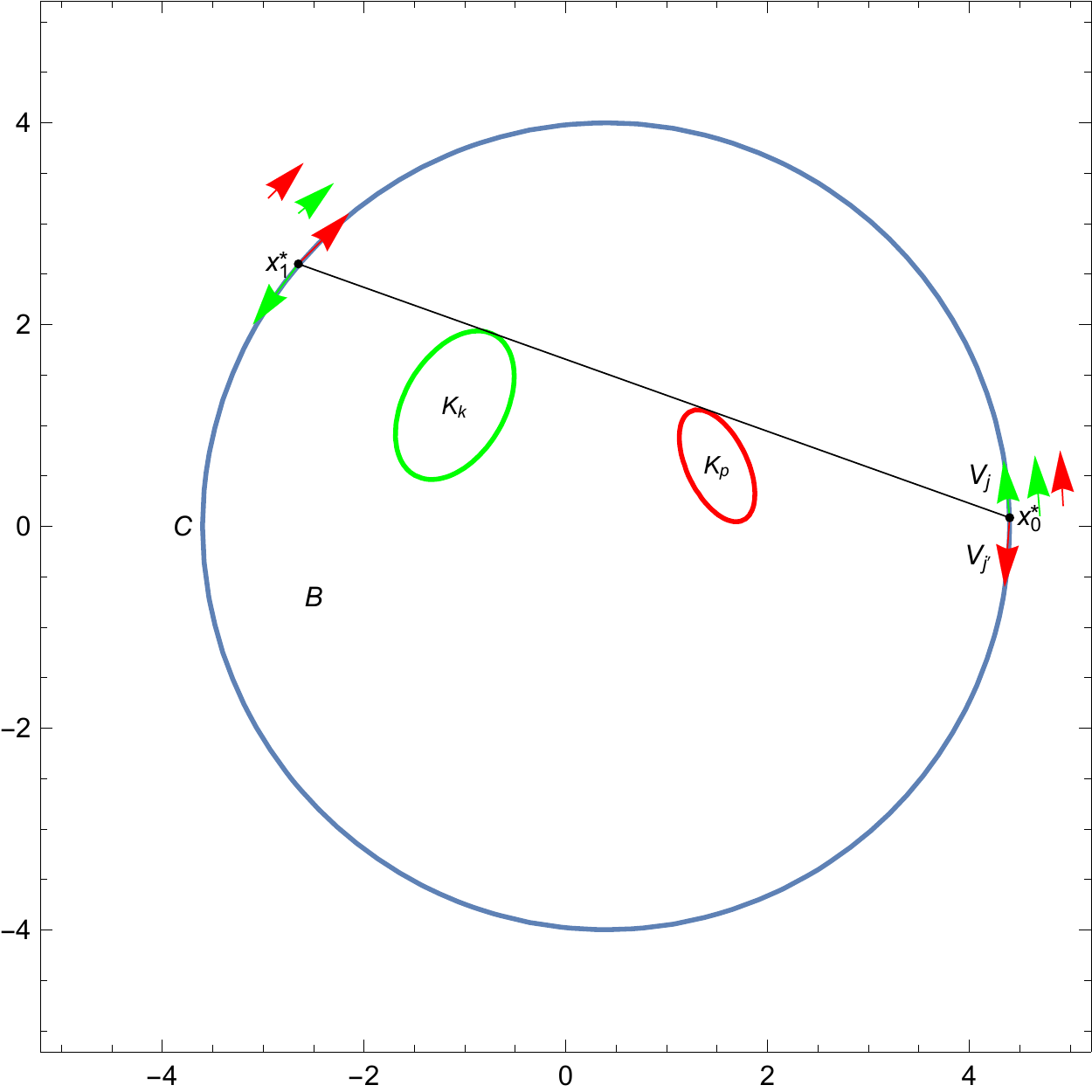} 
   \caption{A linear bitangent (proof of Proposition \ref{prop4})}
   \label{pic3}
\end{figure}

\begin{figure}[h] 
   \centering
   \includegraphics[width=3.1in]{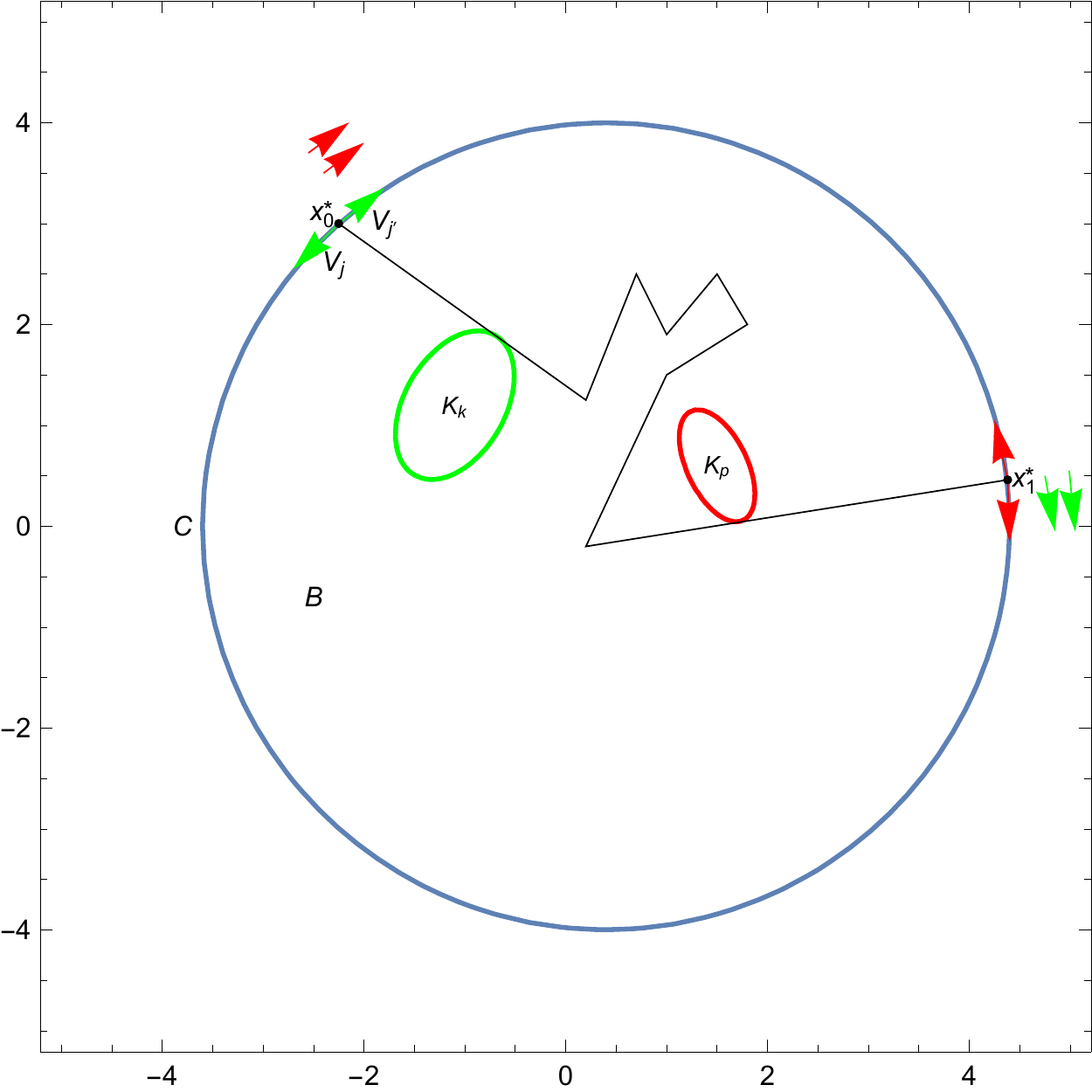} 
   \caption{A nonlinear bitangent geodesic (proof of Proposition \ref{prop4})}
   \label{pic2}
\end{figure}

\spp
When the bitangent  geodesic $x^*$ is linear, there is an open arc $V_{j}\subset C$ of initial points of perturbations  initially tangent to $\partial K_k$, and another open arc $V_{j'}\subset C$ of initial points of perturbations initially tangent  to $\partial K_p$, as in Figure \ref{pic3}, where $x_0^*,V_j,V_{j'}$ appear on the right of the illustration. Perturbations whose initial points are in $V_j$ (green) and $V_{j'}$ (red) have no other points of contact with $\partial K$. 
There are also two unlabelled open arcs $V_{j''},V_{j'''}\subset C$ bordered by $x_0^*$, consisting of initial points of geodesics  whose first points of contact are nontangent to $\partial K$, and whose second points of contact are tangent to $\partial K_k$ (green) or  $\partial K_p$ (red) respectively\footnote{Similarly, the green and red arrows on  the left of Figures \ref{pic3}, \ref{pic2} indicate intervals of terminal points of perturbations.}. 

\spp
Evidently $j\not= j'$, because $V_j$ and $V_{j'}$ are on opposite sides of $x_0^*$, and similarly $j'\not= j'',j'''$ in Figure \ref{pic3}. Indeed,  from the geometry of perturbations of $x^*$, all of $j,j',j'', j'''$ are distinct.  

\spp
In Figure \ref{pic2} the nonlinear bitangent geodesic $x^*$ is tangent to $\partial K_k$ and $\partial K_p$ at the first and last points of contact respectively. It is not tangent anywhere else to $\partial K$, but is reflected at other points of contact, as suggested by the illustration. As before, the nonlinear bitangent is perturbed while maintaining tangency either with $\partial K_k$ (green) or with $\partial K_p$ (red), but now the first and last points of  contact remain on $\partial K_k$ and $\partial K_p$ respectively. The initial points of perturbations tangent to $\partial K_k$ sweep out open arcs $V_j,V_{j'}\subset C$ (green) on either side of $x_0^*$. Initial points of perturbations tangent to $\partial K_p$ give the other intervals $V_{j''},V_{j'''}$ on one side of $x_0^*$, as indicated by the two red arrows on the left of Figure \ref{pic2}. Again $j,j',j'',j'''$ are distinct. 

\spp 
An element $(x_0,u_0,x_1,t)$ of $\tilde \beta _j$ corresponds precisely to the point of tangency (first or last contact) of $x_{x_0, t u_0}$ with $\partial K$. Because there is only one point of tangency it corresponds diffeomorphically to $x_0$. This proves {\it 3}.  

\spp
So $(x_0^*,u_0^*,x_1^*,t^*)$ is an endpoint of precisely $4$ open arcs and {\it 4.} is proved. 

\spp 
Because $\tilde {\mathcal{T}}^1$ is open in $\tilde {\mathcal{T}}^+$, $\cup _{j\geq 1}\partial \tilde \beta _j=\partial \tilde {\mathcal{T}}^1\subseteq \tilde {\mathcal T}^2$. Because $\tilde {\mathcal{T}}^1$ is dense in $\tilde {\mathcal{T}}^+$, $\tilde {\mathcal T}^2=\cup _{j\geq 1}\partial \tilde \beta _j$, proving {\it 5}. 
\qedsymbol
\begin{corollary} $\tilde {\mathcal T}^1$ is the smooth part of the $1$-dimensional space $\tilde {\mathcal T}^+$. \qedsymbol 
\end{corollary}

\begin{remark} For $j\geq 1$ the $o(x_{x_0,tu_0})$ for $(x_0,u_0,x_1,t)\in \tilde \beta _j$ depend only on $\tilde \beta _j$. So we may write them as $o(\tilde \beta _j)$. \qedsymbol
\end{remark} 

\spp
We need the following definitions:
\begin{itemize}
\item The open arcs $\tilde \beta _j\in \tilde {\mathcal B}$ where $1\leq j\leq 4(n^2-n)$ are said to be {\em vacuous}. 
\item For $(x_0,u_0)\in E^2\times S^1$ denote the undirected line through $x_0$ parallel to $u_0$ by $\lambda (x_0,v_0)$. 
\item For $j\geq 1$ define $\lambda _j:V_j\rightarrow \R P^2$ by $\lambda _j(x_0)=\lambda (x_0,u_0)$ where $\psi _j(x_0)=(x_0,u_0,x_1,t)$. 
\item  Denote the envelope of $\lambda _j$ by ${\Lambda}_{\tilde \beta _j}:V_j\rightarrow E^2$. 
\end{itemize}
\section{Extendible Arcs and the Inductive Step}\label{prsec}

\spp
At the end of \S \ref{sec1} the travelling-time data ${\mathcal T}$ is used to find $4(n^2-n)$ open arcs $\beta _j\subset {\mathcal T}_0^1$. 
Each of these is augmented, as described in Proposition \ref{prop4}, to a vacuous open arc $\tilde \beta _j\subset \tilde {\mathcal T}^1$. From the definition in \S \ref{sec1} of the conjugate $\bar j$ of $j$, for $1\leq j\leq 4(n^2-n)$,  
\begin{equation}\label{conjeq}
\Lambda _{\tilde \beta _j}(V_j)~=~\Lambda _{\tilde \beta _{\bar j}}(V_{\bar j}).
\end{equation}
We also obtain $C^\infty$ parameterisations 
$\psi _j:V_j\rightarrow \tilde \beta _j$. More generally (inductively) suppose we have this kind of information where possibly $j>4(n^2-n)$. 

\spp
In precise terms, suppose we are given a $C^\infty$ parameterisation $\psi _j:V_j\rightarrow \tilde \beta _j$ of some possibly nonvacuous arc $\tilde \beta _j\in \tilde {\mathcal B}$. Here $V_j\subset C$ is a maximal open arc with the property that, for all $x_0\in V_j$ and $(x_0,u_0,x_1,t):=\psi _j(x_0)$, the first segment of the geodesic $x_{x_0,tu_0}$ is tangent to $\partial K_k$. The inductive step extends the open arc $\Lambda _{\tilde \beta _j}(V_j)\subset \partial K$ by adjoining another such arc to its clockwise endpoint, as follows. 

\spp
For $x_0^*\in C$ the clockwise terminal limit of $x_0\in V_j$, set 
$$(x_0^*,u_0^*,x_1^*,t^*):=\lim_{x_0\rightarrow x_0^*}\psi _j(x_0)\in \tilde {\mathcal T}^2.$$ 
By Proposition \ref{prop4} there are three other open arcs $\tilde \beta _{j'},\tilde \beta _{j''},\tilde \beta _{j''}\in \tilde {\mathcal B}$  adjacent to $\tilde \beta _j\subset \tilde {\mathcal T}^1$ at $(x_0^*,u_0^*,x_1^*,t^*)$, and the unordered set $\tilde {\mathcal B}_j:=\{ \tilde \beta _{j'},\tilde \beta _{j''},\tilde \beta _{j'''}\}\subset \tilde{\mathcal B}$ is found by inspecting 
$\tilde{\mathcal T}^+$. In the proof of Proposition \ref{prop4}, the arcs $\tilde \beta _j,\tilde \beta _{j'}$ (respectively $\tilde \beta _{j''},\tilde \beta _{j'''}$) are generated by geodesics whose first (respectively last) segments are tangent to $\partial K$. Construct\footnote{From the proof of Proposition \ref{prop4}, $\tilde {\mathcal B}_j^*$ has size $1$ or $3$. } 
$$\tilde {\mathcal B}_j^*~:=~\{ \tilde \beta \in \tilde {\mathcal B}_j: V\cap V_j=\emptyset \}$$
where $V:=\{ x_0: (x_0,u_0,x_1,t)\in \tilde \beta\}$. 
\begin{definition} $\tilde \beta \in \tilde {\mathcal B}_j^*$ is an {\em extension} of $\tilde \beta _j$ 
when the closure $\Lambda _{j,\tilde \beta}$ of $\Lambda _{\tilde \beta _j}(V_j)\cup \Lambda _{\tilde \beta }(V)$ is a $C^\infty$ strictly convex arc in $E^2$. When an extension of $\tilde \beta _j$ exists, the arc $\tilde \beta _j$ is said to be {\em extendible} (otherwise {\em nonextendible}). \qedsymbol
\end{definition}

\begin{proposition}\label{prop5} If $\tilde \beta _j$ is extendible the extension $\tilde \beta \in \tilde {\mathcal B}_j^*$ is unique, and $\Lambda _{\tilde \beta }(V)$ is  an arc in $\partial K_k$. If $\tilde \beta _j$ is nonextendible then $\tilde \beta _j$ is vacuous and $\tilde \beta _{\bar j}$ is extendible. 
\end{proposition} 

\spp
{\bf Proof:} By continuity of $\psi _j$, the bitangent $x_{x_0^*,t^*u_0*}$ is tangent to $\partial K_k$ at some $q:=x_{x_0^*,t^*u_0*}(t_k)$ where $0<t_k<1$, and $q$ is a limit of points of first tangency and first contact with $\partial K_k$. By general position $q$ is either the first point of contact of the bitangent with $\partial K$ or the second point of contact. 

\spp
If $q$ is the first point of contact then $\tilde \beta _j$ is extended by $\tilde \beta _{j'}$ whose associated geodesics maintain tangency to $\partial K_k$. Evidently $\Lambda _{\tilde \beta _{j'}}(V_{j'})$ is an arc in $\partial K_k$. 

\spp
For $j^*=j''$ or $j^*=j'''$, and $(\tilde x_0,\tilde u_0,\tilde x_1,\tilde t)\in \tilde \beta _{j^*}$ near $(x_0^*,u_0^*,x_1^*,t^*)$, the last points of contact of $x_{\tilde x_0,\tilde t\tilde u_0}$ are tangent to $\partial K$ near $q'\in \partial K_{k'}$ where $q'\not= q$. By the argument in \S 3 of \cite{NS1}, the $\lambda (\tilde x_0,\tilde u_0)$ are not all tangent to a $C^\infty$ strictly convex arc, namely $\Lambda _{\tilde \beta _{j^*}}$ is not strictly convex, and $\tilde \beta _{j^*}$ does not extend $\tilde \beta _{j}$. So the extension $\tilde \beta =\tilde \beta _{j'}$ is unique.

\spp
If alternatively  $q$ is the second point of contact, then the first point of tangency is at $q':=x_{x_0^*,t^*u_0*}(s')\in \partial K_{k'}$ where $0<t_{k'}<t_k$ with $k'\not= k$. By Lemma \ref{lem0}, $q'$ is the first point of contact of the bitangent with $\partial K$. By Lemma \ref{lem0}, and because $q$ is the second point of tangency, the bitangent is linear with $q,q'$ the only points of contact with $\partial K$. So $1\leq j\leq 4(n^2-n)$, and $q'$ is the first point of contact of the linear bitangent  $x_{x_1^*,-t^*u_0^*}$ with $\partial K$. Then $\tilde \beta _{\bar j}$ is extended by requiring tangency to $\partial K_{k}$ of the associated geodesics. \qedsymbol

\spp
The arc $\Lambda _{\tilde \beta _j}(V_j)$ in $\partial K_k$ is therefore extended by an incremental arc $\Lambda _{\tilde \beta }(V)$, where $\tilde \beta$ is an extension either of $\tilde \beta _j$ or of $\tilde \beta _{\bar j}$. 
This completes the inductive step. 

\spp
Now the construction of  $\partial K$ proceeds as follows. First $\tilde \beta _j$ is chosen with $1\leq j\leq 4(n^2-n)$, and the inductive step is carried out repeatedly with $\tilde \beta$ replacing $\tilde \beta _j$ after each step, until the incremental arcs $\Lambda _{\tilde \beta }(V)$ in $\partial K$ are acceptably\footnote{Countably many repetitions would be needed for perfect reconstruction.} small. 
Then another vacuous arc is used to restart the iterative process. This is  repeated until all the vacuous arcs are used. Finally $\partial K$ is the union of the closures of all the arcs (initial and incremental) in $\partial K$. 

\spp

\end{document}